\documentclass[12pt]{article}

\usepackage[english]{babel}
\usepackage{amsmath}
\usepackage{amssymb,amsfonts,amsthm}
\usepackage{tikz,color,pgf,graphicx,url,cite}
\usepackage{xspace}

\usetikzlibrary{calc}
\usetikzlibrary{decorations.pathreplacing}
\usetikzlibrary{decorations.pathmorphing}

\usepackage{enumerate}     
\usepackage[T1]{fontenc}
\usepackage{xcolor} 
\usepackage{hyperref}
\hypersetup{
	colorlinks   = true, 
	urlcolor     = black, 
	linkcolor    = black, 
	citecolor   = black 
}

\usepackage[
top=3cm,
bottom=3cm,
inner=2.5cm,      
outer=2.5cm,
footskip=40pt     
]{geometry}

\newtheorem{theorem}{Theorem}[section]
\newtheorem{corollary}[theorem]{Corollary}
\newtheorem{lemma}[theorem]{Lemma}
\newtheorem{proposition}[theorem]{Proposition}

\newtheorem{definition}[theorem]{Definition}

\def\cp{\,\square\,}

\newcommand{\ftmv}{{\sc Fault-tolerant mutual-visibility}\xspace}

\DeclareMathOperator{\diam}{diam}

\DeclareMathOperator{\hull}{hull}
\newcommand{\mv}[1]{\mathrm{f}\mu^{(#1)}}

\newcommand{\vertex}{\node[vertex]}
\tikzstyle{vertex}=[circle, draw, inner sep=0pt, minimum size=6pt]

\title{Fault-tolerant mutual-visibility: complexity and solutions for grid-like networks}

\author{
Serafino Cicerone $^{a,}$\thanks{Email: \texttt{serafino.cicerone@univaq.it}}
\and
Gabriele {Di Stefano} $^{a,}$\thanks{Email: \texttt{gabriele.distefano@univaq.it}}
\and
Sandi Klav\v zar $^{b,c,d,}$\thanks{Email: \texttt{sandi.klavzar@fmf.uni-lj.si}}
\and Gang Zhang $^{c,e,}$\thanks{Email: \texttt{gzhang@stu.xmu.edu.cn}}
}

\begin{document}

\maketitle
 
\begin{center}
	$^a$ Department of Information Engineering, Computer Science, and Mathematics,
	     University of L'Aquila, Italy \\
	\medskip

	$^b$ Faculty of Mathematics and Physics, University of Ljubljana, Slovenia\\
	\medskip
	
	$^c$ Institute of Mathematics, Physics and Mechanics, Ljubljana, Slovenia\\
	\medskip
	
	$^d$ Faculty of Natural Sciences and Mathematics, University of Maribor, Slovenia\\
	\medskip
	
	$^e$ School of Mathematical Sciences, Xiamen University, China
\end{center}

\begin{abstract}
Networks are often modeled using graphs, and within this setting we introduce the notion of $k$-fault-tolerant mutual visibility. Informally, a set of vertices $X \subseteq V(G)$ in a graph $G$ is a $k$-fault-tolerant mutual-visibility set ($k$-ftmv set) if any two vertices in $X$ are connected by a bundle of $k+1$ shortest paths such that:
($i$) each shortest path contains no other vertex of $X$, and
($ii$) these paths are internally disjoint.
The cardinality of a largest $k$-ftmv set is denoted by $\mathrm{f}\mu^{k}(G)$. The classical notion of mutual visibility corresponds to the case $k = 0$.

This generalized concept is motivated by applications in communication networks, where agents located at vertices must communicate both efficiently (i.e., via shortest paths) and confidentially (i.e., without messages passing through the location of any other agent). The original notion of mutual visibility may fail in unreliable networks, where vertices or links can become unavailable.

Several properties of $k$-ftmv sets are established, including a natural relationship between $\mathrm{f}\mu^{k}(G)$ and $\omega(G)$, as well as a characterization of graphs for which $\mathrm{f}\mu^{k}(G)$ is large. It is shown that computing $\mathrm{f}\mu^{k}(G)$ is NP-hard for any positive integer $k$, whether $k$ is fixed or not. Exact formulae for $\mathrm{f}\mu^{k}(G)$ are derived for several specific graph topologies, including grid-like networks such as cylinders and tori, and for diameter-two networks defined by Hamming graphs and by the direct product of complete graphs.
\end{abstract}


\medskip\noindent
\textbf{AMS Math.\ Subj.\ Class.\ (2020)}: 05C12, 05C69, 05C76, 68Q17


\section{Introduction}
%
Graphs provide a powerful framework for modeling real-world networks, from communication and transportation systems to social and biological structures. In many applications, such as the coordination of autonomous agents or the deployment of sensors in a network, vertices represent locations, and the ``visibility'' between them is a critical operational constraint. A recently introduced concept, \textbf{mutual-visibility}, provides a formal graph-theoretic model for this idea \cite{distefano-2022}.

Formally, given a connected graph $G$ and a subset of vertices $X\subseteq V(G)$, two vertices $u, v \in V(G)$ are said to be \emph{visible} to each other with respect to $X$ if there exists a shortest $u,v$-path in $G$ that contains no internal vertices from $X$. A set $X \subseteq V(G)$ is then defined as a \emph{mutual-visibility set} if every pair of distinct vertices in $X$ is visible with respect to $X$. The \emph{mutual-visibility number}, $\mu(G)$, is the cardinality of a largest mutual-visibility set in $G$. This concept effectively models a scenario where agents, located at the vertices of $S$, can ``see'' each other if a line-of-sight (i.e., a shortest path) is not ``blocked'' by another agent.

This definition of mutual-visibility, however, implicitly assumes a perfectly reliable network. The visibility between two agents $u$ and $v$ depends on the existence of \emph{at least one} ``clear'' shortest path, where clear means that these paths do not pass through any \emph{other} agents located in $X$. In practical applications, networks are often unreliable since vertices or links may fail. If the single clear path between $u$ and $v$ is compromised by a faulty vertex (that is not in $X$), the agents lose visibility, even if alternative paths exist. This vulnerability motivates the need for a more robust, fault-tolerant variant of mutual-visibility.

The principle of using path redundancy to achieve fault tolerance is a cornerstone of graph theory and network design. This concept is fundamentally rooted in Menger's theorem, which equates the $k$-connectivity of a graph to the existence of $k$ internally disjoint paths between any two vertices, cf.~\cite{beineke-2012}. This principle has been extensively applied to model and solve problems in reliable systems. For example, in \emph{survivable network design}, the goal is to create network topologies that maintain connectivity properties despite component failures, often by ensuring a minimum number of disjoint paths \cite{Kerivin-2005}. Similarly, \emph{fault-tolerant routing} protocols are designed to find alternative paths to bypass faulty nodes or links, guaranteeing message delivery \cite{Gomez-2006}. These fields demonstrate a clear consensus: resilience against $k$ failures can be achieved by providing $k+1$ disjoint resources.

In this paper, we apply this established principle of fault tolerance to the mutual-visibility framework. We introduce the \textbf{$k$-fault-tolerant mutual-visibility}, a concept that strengthens the original definition by requiring a specified redundancy in the visibility paths.
Informally, a set of vertices $X$ is a $k$-fault-tolerant mutual-visibility set (or $k$-ftmv set) if any two agents (vertices) within $X$ are connected by a bundle of $k+1$ clear shortest paths. These paths must be \emph{internally disjoint}, meaning they share no internal vertices. Moreover, as in the original definition, clear means that these paths do not pass through any \emph{other} agents located in $X$. This built-in redundancy ensures that the mutual-visibility between any two agents is robust: it can withstand the failure of up to $k$ internal vertices and still be maintained by at least one remaining clear shortest path.

The cardinality of a largest $k$-ftmv set of a graph $G$ is denoted by $\mv{k}(G)$. Since $\mv{0}(G) = \mu(G)$ by definition, it follows that the newly introduced concept captures and extends the standard notion of mutual-visibility.

\paragraph{Related works.} 
Questions regarding mutual-visibility and related arrangements of points in the Euclidean plane have been a subject of study since the late 19th century, highlighted by classical problems like the no-three-in-line problem~\cite{Hardy-2008}. As a counterpart of the classic notion of visibility in the Euclidean plane, the graph theoretic notion of mutual-visibility has been introduced in~\cite{distefano-2022}. 
Despite being introduced recently, the mutual-visibility in graphs gathered significant interest within the research community becoming a flourishing line of research (e.g., see~\cite{%
axenovich-2025+, 
boruzanli-2024,
BresarYero-2024,
Bujtas-2025,
cicerone-2025,
cicerone-2025a,
CiDiDrHeKlYe-2023,
cicerone-2023,
cicerone-2024b,
korze-2024,
korze-2024+,
klavzar-2025,
kuziak-2023,
roy-2025,
roy-2025b,
tian-2024}). 
Since in~\cite{distefano-2022} it is shown that computing $\mu(G)$ is a NP-complete problem, in most of such papers, the main aim was solving the mutual‑visibility problem for various classes of graphs. It is a remarkable fact that, from a combinatorial viewpoint, the mutual‑visibility problems connects several topics: in Cartesian products of complete graphs, it is equivalent to instances of Zarankiewicz’s problem~\cite{cicerone-2023}; bounds on smallest maximal mutual-visibility sets relate to the Bollob\'as–Wessel theorem~\cite{BresarYero-2024}; and the various optimization problems arising from mutual-visibility variants can be reformulated as Turán‑type problems on hypergraphs and line graphs~\cite{Bujtas-2025,cicerone-2024b}.

Other directions of research have focused on possible extensions of the definition of mutual-visibility. In~\cite{CiDiDrHeKlYe-2023}, variations of the initial definition are introduced, based on the extension of the visibility property of vertices that are in and/or outside $X$. Such variations are called \emph{total}, \emph{outer} and \emph{dual} mutual-visibility problems. 
In~\cite{bresar-2025,klavzar-2025}, the authors investigate proper vertex colorings of graphs whose color classes are mutual-visibility sets. In~\cite{kuziak-2025}, the concept of a \emph{$d$-visibility set} is introduced---a set $X$ in which each pair of vertices either has distance greater than $d$, or is connected by a shortest path whose internal vertices are not in the set. In this way, the original mutual-visibility concept is restricted to vertices that are close to each other within the given threshold $d$.

The original mutual-visibility definition is also naturally related to that of \emph{general position}. A set of vertices $X$ is a general position set if no three vertices from $X$ lie in a common shortest path~\cite{manuel-2018,manuel-2018a,parthasarathy-2016}. The distinction between a general‑position set and a mutual‑visibility set lies in the nature of the restriction applied to shortest paths. While mutual-visibility calls for the existence of at least one clear shortest-path between them, the general‑position set, by contrast, imposes a stronger restriction asking that all the shortest paths between them are clear.

\paragraph{Results.}
We first formally introduce the new concept of fault-tolerant mutual-visibility in graphs and provide some basic properties. Among them, we show a natural relationship between $\mv{k}(G)$ and $\omega(G)$. Concerning computational issues, we formally prove that computing $\mv{k}(G)$ is a NP-hard problem for any fixed and not fixed positive integer $k$.

Then, we study the structure of the graphs whose fault-tolerant mutual-visibility number is large.  More precisely, we describe the structure of the graphs $G$ with $\mv{k}(G) = n(G) - t$ for all $k\ge 1$ and $t\in \{0,1,\dots, k+1\}$, showing that the graphs with $\mv{k}(G) = n(G) - t$ coincide with the graphs whose clique number is $n(G) - t$. 

Concerning computing the fault-tolerant mutual-visibility number of relevant graph families, we show exact formulae for $\mv{k}(G)$ when $G$ corresponds to a grid, a cylinder, or a torus graph. We also determine the fault-tolerant mutual-visibility number of some classes of graphs of diameter two, like as Hamming graphs and the direct product of complete graphs.

\section{Basic terminology and notation}\label{sec:notation}
%
All graphs considered in this paper are finite and simple. They are also connected unless stated otherwise. Given a graph $G = (V(G), E(G))$, its order will be denoted by $n(G)$. If $v\in V(G)$, then $N_G(v)$ denotes the set of {\em neighbors} of $v$,  that is, the set of vertices that are adjacent to $v$; $N_G[v] = N_G(v)\cup \{v\}$ contains the {\em closed neighbors}. If $X\subseteq V(G)$, then $N_G(X)=\bigcup_{v\in X} N(v)$. The \emph{degree} of $v$ is $\deg_G(v)=|N_G(v)|$ and $\Delta(G)$ denotes the maximum degree of $G$. A \emph{universal vertex} is a vertex that is adjacent to all other vertices of the graph, that is $\deg_G(v)=n(G)-1$. If $X\subseteq V(G)$, then $G[X]$ denotes the subgraph of $G$ induced by $X$, that is the maximal subgraph of $G$ with vertex set $X$. The subgraph of $G$ induced by $V(G)\setminus X$ is denoted by $G-X$, and by $G-v$ when $X=\{v\}$. A \emph{cut-vertex} of a graph $G$ is a vertex whose removal increases the number of connected components. A \emph{block} of $G$ is a maximal connected subgraph of $G$ that has no cut-vertex.

The distance function $d_G$ on a graph $G$ is the usual shortest-path distance. The longest distance in $G$ is its \emph{diameter}, denoted by $\diam(G)$. The subgraph $G'$ of $G$ is \emph{convex} if, for every two vertices of $G'$, every shortest path in $G$ between them lies completely in $G'$.  The \emph{convex hull} of $V'\subseteq V(G)$, denoted as $\hull(V')$, is defined as the smallest convex subgraph containing $V'$. 

We use standard notation for classic graphs with $n$ vertices: $P_n$ is the path graph, $C_n$ is the cycle graph, $K_n$ is the complete graph. Moreover, $K_{m,n}$ is the complete bipartite graph of order $n+m$, and, as a special case, $K_{1,n}$ is the star with $n$ pendant vertices.

A \emph{clique} is an induced subgraph in which every two distinct vertices are adjacent. The \emph{clique number} $\omega(G)$ is the number of vertices in a maximum clique in $G$. An \emph{independent set} is a set of vertices in a graph, no two of which are adjacent. The \emph{independence number} $\alpha(G)$ is the size of the largest possible independent set in $G$.

Given two graphs $G$ and $H$, their {\em Cartesian product} and {\em direct product} are respectively denoted by $G\cp H$ and $G\times H$. Each of these products has the vertex set $V(G)\times V(H)$. Vertices $(g,h)$ and $(g',h')$ are adjacent in $G\cp H$ if either $g=g'$ and $hh'\in E(H)$, or $gg'\in E(G)$ and $h=h'$. Vertices $(g,h)$ and $(g',h')$ are adjacent in $G\times H$ if $gg'\in E(G)$ and $hh'\in E(H)$. If $h\in V(H)$ and $*\in \{\cp, \times\}$, then the subgraph of $G * H$ induced by the vertex set $\{(g,h):\ g\in V(G)\}$ is called an {\em $G$-layer} of $G*H$ and denoted by $G^h$. For $g\in V(G)$, the {\em $H$-layer} $^gH$ is defined analogously. 

\section{Fault-tolerant mutual-visibility and its general properties}
\label{sec:concept+properties}

In this section, we formally introduce the key concepts of this paper and deduce some of the basic properties. Let's get started right away with the key definitions.  

\begin{definition}\label{def:fault-tolerant-mv}
Let $k\ge 0$. Then a set $X\subseteq V(G)$ is a {\em $k$-fault-tolerant mutual-visibility set} (shortly {\em $k$-ftmv set}) if for any non-adjacent $u,v\in X$ there exist $k+1$ internally disjoint shortest $u,v$-paths $Q_i$, $i\in [k+1]$, such that $V(Q_i)\cap X = \{u,v\}$. The cardinality of a largest $k$-ftmv set of $G$ is denoted by $\mv{k}(G)$, and a largest $k$-ftmv set is referred as a {\em $\mv{k}(G)$-set}. 
\end{definition}

A set $X$ is thus a $k$-ftmv set if, whenever $k$ vertices from $V(G)\setminus X$ are removed, for each non-adjacent $u,v\in X$ there remains a path $Q$ of length $d_G(u,v)$ such that $V(Q)\cap X = \{u,v\}$. Informally, the mutual-visibility of the vertices from $X$ is not affected by having at most $k$ faulty vertices. 

Observe that $0$-ftmv sets coincide with the standard mutual-visibility sets, in particular, $\mv{0}(G) = \mu(G)$. Note also that $\mv{k}(G)$ is well-defined for any $k\ge 0$ because the set of vertices of an arbitrary clique of $G$ forms a $k$-ftmv set of $G$ for any $k$. In particular, 
\begin{equation}
\label{eq:omega-is-lower-bound}   
\mv{k}(G)\ge \omega(G),\ k\ge 0. 
\end{equation}
From the definition we also get, 
\begin{equation}
\label{eq:monotone}
\mv{k}(G) \le \mv{k-1}(G),\ k\ge 1. 
\end{equation}

We continue with the following two general properties. 

\begin{proposition}\label{prop:delta}
Let $G$ be a graph. If $k\ge\Delta(G)$ then  $\mv{k}(G)= \omega(G)$.
\end{proposition}

\begin{proof}
Let $X$ be a $k$-ftmv set of $G$. By contradiction, assume there exist two non-adjacent vertices $u,v\in X$.  There are at most $\deg_G(u)$ internally disjoint (shortest) $u,v$-paths in $G$ and since $\deg_G(u) \le \Delta(G)\le k$,  they are at most $k$. Hence, all the vertices in $X$ are adjacent and the result follows.
\end{proof}

\begin{proposition}
\label{prop:cut-vertex}
Let $k\ge 1$, and let $B_1,\dots, B_\ell$ be the blocks of a nontrivial connected graph $G$. If $X$ is a $\mv{k}(G)$-set of $G$, then $X\subseteq V(B_i)$ for some $i\in [\ell]$. In particular, $$\mv{k}(G) = \max\ \{\mv{k} (B_i):\ i \in [\ell]\}\,.$$
\end{proposition}

\begin{proof}
Let $X$ be a $\mv{k}(G)$-set. Suppose that two vertices $u$ and $v$ of $X$ do not belong to the same block of $G$. Let $u\in X\cap B_i$ and $v\in X\cap B_j$, and let $Q$ be a shortest $u,v$-path in $G$. Since $i\ne j$, the path $Q$ contains at least one cut-vertex $w$, where $w\ne u,v$. Then $w\notin X$ and therefore every shortest $u,v$-path passes $w$. As this is not possible, we conclude that $X$ must be contained in a single block of $G$. Hence in particular, $\mv{k}(G)=|X| = \max_{i \in [\ell]}\{\mv{k}(B_i)\}$.     
\end{proof}

Notice that, for each non-trivial tree $T$, Proposition~\ref{prop:cut-vertex} implies $\mv{k}(T)=2$ for each $k\ge 1$. Conversely, in~\cite[Corollary 4.3]{distefano-2022} it is shown that $\mv{0}(T)=\mu(T)=|L|$, where $L$ is the set of leaves of $T$. It is worth also remarking that in~\cite[Lemma 2.5]{distefano-2022} it is shown that a $\mv{0}$-set of $G$ containing no cut-vertices always exists.

The following lemmas were formulated in~\cite{distefano-2022} for the original version of mutual-visibility. We show that they can be extended to the $k$-fault-tolerance cases, $k\ge 1$.

\begin{lemma}\label{lem:convex-subgraph}
Let $G$ be a graph and $H$ a convex subgraph of $G$. Then $\mv{k}(H) \leq \mv{k}(G)$. Moreover, if $X$ is a $k$-ftmv set of $G$, then $X\cap V(H)$ is a $k$-ftmv set of $H$.
\end{lemma}
\begin{proof}
Since every $k$-ftmv set of $H$ is a $k$-ftmv set of $G$, the first statement holds.

Let $u$ and $v$ be two non-adjacent vertices in $X' = X\cap V(H)$. Since $X$ is a $k$-ftmv set of $G$, there are $k+1$ internally disjoint shortest $u,v$-paths $Q_i$, $i\in [k+1]$, such that $V(Q_i)\cap X = \{u,v\}$. By the definition of convex subgraph, all such $k+1$ shortest $u,v$-paths are in $H$. Since they have no internal vertices in $X$, they also have no internal vertices in $X'$. This implies that $X'$ is a $k$-ftmv set of $H$.
\end{proof}

\begin{lemma}\label{lem:covering-by-convex-subgraphs}
Let $G$ be a graph and $V_1,V_2,\ldots,V_{\ell}$ subsets of $V(G)$ such that $\bigcup_{i=1}^{\ell} V_i = V(G)$. Then, $\mv{k}(G) \le \sum_{i=1}^{\ell} \mv{k}( \hull(V_i) ) $, for each $k\ge 0$.
\end{lemma}
\begin{proof}
Assume $\mv{k}(G) > \sum_{i=1}^{\ell} \mv{k}( \hull(V_i) ) $ and let $X$ be a $\mv{k}(G)$-set. Since $\bigcup_{i=1}^{\ell} V_i = V(G)$, each elements of $X$ is in at least one set $\hull(V_i)$. For each $i\in [\ell]$, let $X_i$ be the set of vertices that are in $X$ and in $\hull(V_i)$. Then $\sum_{i=1}^{\ell} |X_i| \ge |X| = \mv{k}(G) > \sum_{i=1}^{\ell} \mv{k}( \hull(V_i) )$. Hence, there must exist a set $X_j$, $j\in [\ell]$, such that $|X_j| > \mv{k}(\hull(V_j))$. This is a contradiction because Lemma~\ref{lem:convex-subgraph} implies that $X_j$ is a $k$-ftmv set of $\hull(V_j)$ and its size cannot be larger than $\mv{k}(\hull(V_j))$.
\end{proof}

Let $G$ be a connected graph. Then $S\subseteq V(G)$ is a {\em clique cut-set} if $G[S]$ is a clique and $G-S$ is disconnected. A clique cut-set $S$ is {\it minimal} if $S$ does not contain a proper clique cut-set. Using the above two lemmas, the following can be proved. 

\begin{proposition}
\label{propo-2-cuts}
Let $S$ be a minimal clique cut-set with $|S|=2$.  If $G_1,\dots,G_{\ell}$, $\ell \ge 2$, are the components of $G-S$, and $\widehat{G_i} = G[V(G_i)\cup S]$, $i\in [\ell]$, then the following hold. 
\begin{enumerate}
\item[(i)] $\displaystyle{ \max\{\mv{1}(\widehat{G_i}):\ i \in [\ell]\} \le \mv{1}(G) \le \sum_{i=1}^{\ell}\mv{1}(\widehat{G_i})}$, and the bounds are sharp.
\item[(ii)] If $k\ge 2$, then $\mv{k}(G) = \max\{\mv{k}(\widehat{G_i}): i \in [\ell]\}$.
\end{enumerate}
\end{proposition}

\begin{proof}
(i) Since every $\widehat{G_i}$, $i \in [\ell]$, is a convex subgraph of $G$, the lower bound follows by Lemma~\ref{lem:convex-subgraph}, while the upper bound is implied by Lemma~\ref{lem:covering-by-convex-subgraphs}.

To demonstrate the sharpness of the lower bound, consider the graphs $C_{m,n}$, $m, n\ge 4$, obtained from disjoint cycles $C_m$ and $C_n$ by identifying an arbitrary edge of $C_m$ with an arbitrary edge of $C_n$. Then we can verify that $\mv{1}(C_m) = \mv{1}(C_n) = \mv{1}(C_{m,n}) = 2$. Similarly, if $G_{m,n}$, $m\ge 3$, $n\ge 4$, is obtained from $K_m$ and $C_n$ by identifying an arbitrary edge of $K_m$ with an arbitrary edge of $C_n$, then $\mv{1}(C_n) = 2$ and $\mv{1}(K_m) = \mv{1}(G_{m,n}) = m$.

For the sharpness of the upper bound, consider the graphs $H_{m,\ell}$, $m\ge 4, \ell\ge 2$, constructed as follows. First take the disjoint union of a $P_2$, where $V(P_2) = \{x, y\}$, and $\ell$ disjoint copies of $P_m$. Then add all possible edges between $x$ and $y$, and the vertices from the $\ell$ copies of $P_m$. The set $S = \{x, y\}$ is a minimal clique cut-set, and $G - S$ consists of $\ell$ components $G_i$ isomorphic to $P_m$. See Figure~\ref{fig:H-ml} where $H_{5,4}$ is drawn. We can infer that $\mv{1}(\widehat{G_i}) = m$, $i\in [\ell]$, and that $\mv{1}(H_{m,\ell}) = \ell m$.

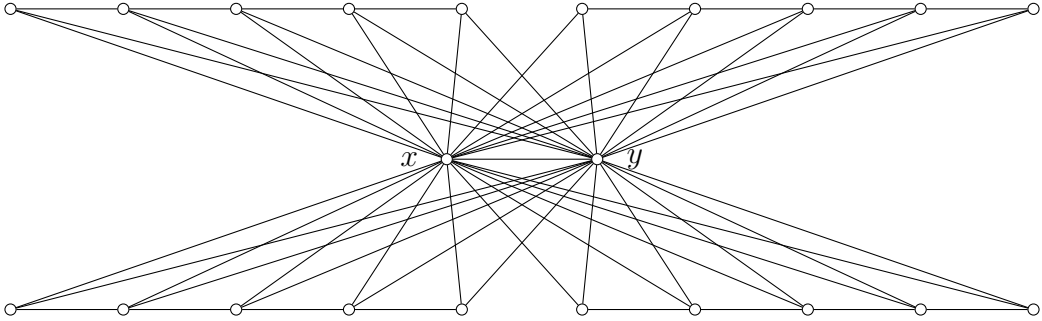
\begin{figure}[ht!]
	\begin{center}
		\begin{tikzpicture}[scale=1.0]
			\tikzstyle{vertex}=[circle, draw, inner sep=0pt, minimum size=6pt]
			\tikzset{vertexStyle/.append style={rectangle}}

                \vertex (1) at (0,0) [scale=0.7,fill=white] {};
                \vertex (2) at (2,0) [scale=0.7,fill=white] {};

                \node (x) at (-0.5,0) {$x$};
                \node (y) at (2.5,0) {$y$};

                \vertex (3) at (-2.8,2) [scale=0.7,fill=white] {};
                \vertex (4) at (-1.3,2) [scale=0.7,fill=white] {};
                \vertex (5) at (0.2,2) [scale=0.7,fill=white] {};
                \vertex (-4) at (-4.3,2) [scale=0.7,fill=white] {};
                \vertex (-5) at (-5.8,2) [scale=0.7,fill=white] {};

                \vertex (6) at (1.8,2) [scale=0.7,fill=white] {};
                \vertex (7) at (3.3,2) [scale=0.7,fill=white] {};
                \vertex (8) at (4.8,2) [scale=0.7,fill=white] {};
                \vertex (-7) at (6.3,2) [scale=0.7,fill=white] {};
                \vertex (-8) at (7.8,2) [scale=0.7,fill=white] {};

                \vertex (9) at (-2.8,-2)[scale=0.7,fill=white] {};
                \vertex (10) at (-1.3,-2) [scale=0.7,fill=white] {};
                \vertex (11) at (0.2,-2) [scale=0.7,fill=white] {};
                \vertex (-10) at (-4.3,-2) [scale=0.7,fill=white] {};
                \vertex (-11) at (-5.8,-2) [scale=0.7,fill=white] {};

                \vertex (12) at (1.8,-2)[scale=0.7,fill=white] {};
                \vertex (13) at (3.3,-2) [scale=0.7,fill=white] {};
                \vertex (14) at (4.8,-2) [scale=0.7,fill=white] {};
                \vertex (15) at (6.3,-2) [scale=0.7,fill=white] {};
                \vertex (16) at (7.8,-2) [scale=0.7,fill=white] {};
                
                \path
                (1) edge (2)
                (1) edge (3)
                (1) edge (4)
                (1) edge (5)
                (1) edge (-4)
                (1) edge (-5)
                (1) edge (6)
                (1) edge (7)
                (1) edge (8)
                (1) edge (-7)
                (1) edge (-8)
                (1) edge (9)
                (1) edge (10)
                (1) edge (11)
                (1) edge (-10)
                (1) edge (-11)
                (1) edge (12)
                (1) edge (13)
                (1) edge (14)
                (1) edge (15)
                (1) edge (16)
                ;

                \path
                (2) edge (3)
                (2) edge (4)
                (2) edge (5)
                (2) edge (-4)
                (2) edge (-5)
                (2) edge (6)
                (2) edge (7)
                (2) edge (8)
                (2) edge (-7)
                (2) edge (-8)
                (2) edge (9)
                (2) edge (10)
                (2) edge (11)
                (2) edge (-10)
                (2) edge (-11)
                (2) edge (12)
                (2) edge (13)
                (2) edge (14)
                (2) edge (15)
                (2) edge (16)
                ;

                \path
                (3) edge (4)
                (4) edge (5)
                (-5) edge (-4)
                (-4) edge (3)
                (6) edge (7)
                (7) edge (8)
                (8) edge (-7)
                (-7) edge (-8)
                (9) edge (10)
                (10) edge (11)
                (9) edge (-10)
                (-10) edge (-11)
                (12) edge (13)
                (13) edge (14)
                (14) edge (15)
                (15) edge (16)
                ;

	    \end{tikzpicture}
		\caption{The graph $H_{5,4}$.}
		\label{fig:H-ml}
	\end{center}
\end{figure}

(ii) Let $k\ge 2$ and let $S$ be such that $|S|=2$. Let further $X$ be a $\mv{k}$-set of $G$. If $u,v\in X$, then since $S$ is a (minimal clique) cut-set, there exists $i\in [\ell]$ such that $u,v\in \widehat{G_i}$. It follows that $X\subseteq V(\widehat{G_i})$. Lemma~\ref{lem:convex-subgraph} completes the argument. 
\end{proof}

A \emph{module} of a (connected) graph $G$ is a set $S\subseteq V(G)$ such that all vertices in $S$ have the same set of neighbors among vertices not in $S$. This is a fundamental concept in graph theory, as evidenced by its many other names, including autonomous sets, homogeneous sets, and  externally related sets, see~\cite{brand-1999}. We say that $S\subseteq V(G)$ is a {\em cut-module} if $G[S]$ is a module and $G-S$ is disconnected. A cut-module $S$ is {\it minimal} if $S$ does not contain a proper cut-module. 

\begin{lemma}\label{lem:cuts}
Let $S$ be a minimal cut-module of a graph $G$ such that $|S|=k+1$. If  $N(S)\cap H > k$ for each  connected component $H$ of $G-S$, then there exists a $k$-ftmv set $X$ of $G$ such that $X\cap S=\emptyset$.
\end{lemma}

\begin{proof}
Let $S$ be a minimal cut-module $S$ of  $G$ such that $|S|=k+1$ and let $X$ be an arbitrary $k$-ftmv set of $G$. There is nothing to prove if $X\cap S=\emptyset$, hence assume that $X\cap S\not =\emptyset$. Then all  vertices in $X\setminus S$ must belong to a single connected component of $G-S$. Let $H'$ be one of the other components such that $X\cap H'=\emptyset$. Then $(X\setminus S)\cup (N(S)\cap H')$ is a fault-tolerant mutual-visibility set of size at least $|X|$ and with no vertex from $S$ as required. 
\end{proof}

\section{Complexity of the related problems}
\label{sec:complexity}

In this section, we assess the computational complexity of computing $\mv{k}(G)$. To this end, we first introduce the corresponding decision problem: given a graph $G$ and a positive integer $t$, the $k$-\ftmv problem (one problem for each possible fixed $k$) asks if $\mv{k}(G)\geq t$.  The version of the problem in which $k$ is part of the input defines the \ftmv problem: given a graph $G$, an integer $k>0$ and a positive integer $t$, it asks if $\mv{k}(G)\geq t$.

\begin{theorem}\label{thm_np-complete}
Given a positive integer $k$, the $k$-\ftmv problem is NP-complete, even for graphs of diameter four.
\end{theorem}

\begin{proof}
Given a subset of vertices of any graph, it is possible to test in polynomial time whether it is a $k$-fault-tolerant mutual-visibility set or not. Consequently, the problem is in NP. We prove that the {\sc Independent Set} problem (equivalent to the {\sc Clique} problem in the complement graph and shown as NP-complete in~\cite{Karp72}), polynomially reduces to the $k$-\ftmv problem. The {\sc Independent Set} problem asks if the independence number $\alpha(G)$ of a graph $G$ is at least a given integer $t$. 

Let $(G,t)$ be an arbitrary instance of the {\sc Independent Set} problem. We assume $G$ not complete, as these assumptions do not affect the hardness of the problem. Since $G$ is connected and not complete, we get $t\geq 2$. Starting from $(G,t)$, we construct an instance $(G',t')$ of the $k$-\ftmv problem as follows.

\begin{figure}[ht!]
\begin{center}
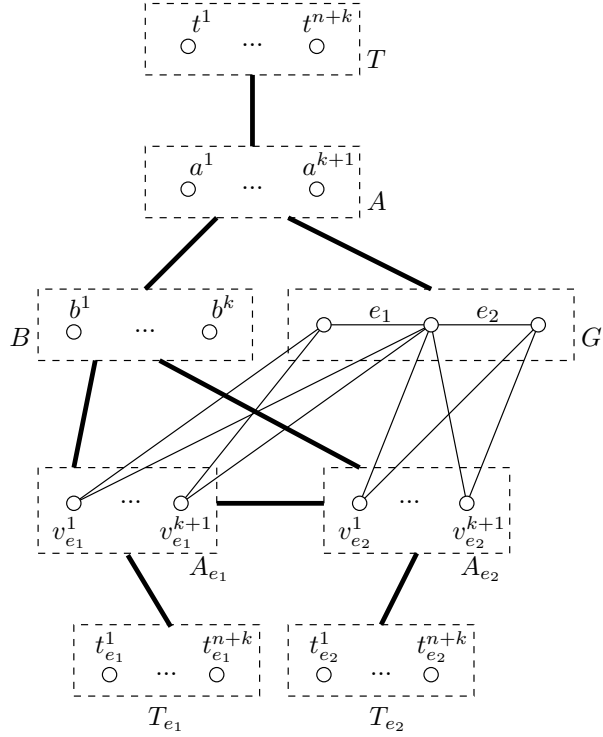
\end{center}
\caption{Visualization of the construction of $G'$, as described in the proof of Theorem~\ref{thm_np-complete}, starting from a graph $G \cong P_3$. Dashed rectangles highlight the subsets of vertices used in the construction. Thin lines between vertices represent real edges. A thick black line between two generic rectangles $R_1$ and $R_2$ simply indicates that every vertex in $R_1$ is adjacent to every vertex in $R_2$.}
\label{fig:np-complete}
\end{figure}

Let $V(G) = [n]$ and $|E(G)|=m$. We add a set of $k+1$ vertices $A=\{a^l:\ l\in [k+1]\}$ each connected to all vertices in $V(G)$. For each edge $e=ij$ of $G$, we add a set of $k+1$ vertices $A_e =\{v_e^l:\ l \in [k+1]\}$  and  edges $iv_e^l$ and $jv_e^l$, for each $l \in [k+1]$. Also, for each pair of edges $e,e'\in E(G)$ we add all edges between $v_e^l$ and $v_{e'}^{l'}$, for each $l,l'\in [k+1]$.  We also add a set of $k$ vertices $B=\{b^l:\ l\in [k]\}$ each connected to all vertices in $A\cup\bigcup_{e\in E(G)} A_e$. Finally, we add a set $T$ of $n+k$ vertices, each connected to all vertices in $A$, and for each edge $e$, we add a set $T_e$ of $n+k$ vertices, each connected to all vertices in $A_e$. This concludes the construction of $G'$. Finally, set $t'= (m+1)(n+k)+t$. Notice that the order of $G'$ is polynomial with respect to the input parameters $n$, $m$, and $t$ (we remind that $k$ is fixed). An example of $G'$ is given in  Figure~\ref{fig:np-complete}. 

In what follows, we prove that $\alpha(G)\geq t$ if and only if $\mv{k}(G')\geq t'$. 

Let $I\subseteq V(G)$ be an independent set of $G$ of size greater than or equal to $t$, and consider the set of vertices $X=I\cup T\cup \bigcup_{e\in E(G)}T_e$ in $G'$. Since $|I|\ge t$, then $|X|\ge (m+1)(n+k)+t = t'$. We now prove that $X$ is a $k$-ftmv set of $G'$. Let $u,v$ be two vertices in $T$. There are $k+1$ shortest $u,v$-paths of length two passing through vertices in $A$. The same holds if $u\in T$ and $v\in I$, or if $u,v\in I$. Let $u\in T$ and $v\in T_e$, $e=ij\in E(G)$. Each path connecting them is of the form $u,x,y,z,v$, where $x\in A$, $y\in B\cup\{i,j\}$, and $z\in A_e$. Since $i$ and $j$ cannot be both in $X$, there are $k+1$ possible choices of $x$, $y$, and $z$, leading to $k+1$ shortest $u,v$-paths of length four without internal vertices in $X$. Let $u\in I$ and $v\in T_e$. If $e$ is incident to $u$, then there are $k+1$ shortest $u,v$-paths of length two passing through vertices in $A_e$. Indeed, for otherwise, if $u$ is a vertex of an edge $e'\not = e$, there are $k+1$ shortest $u,v$-paths of length three of the form $u,x,y,v$, where $x\in A_e$ and $y\in A_{e'}$. Let $u\in T_e$ and $v \in T_{e'}$. If $e=e'$, then there are $k+1$ $u,v$-shortest paths of length two passing through vertices in $A_e$, for otherwise there are $k+1$ $u,v$-shortest paths of length three of the form $u,x,y,v$, where $x\in A_e$ and $y\in A_{e'}$.

In conclusion, each pair of vertices in $X$ is connected by $k+1$ shortest paths without internal vertices in $X$.

\medskip
Let now $X$ be a $k$-ftmv set of $G'$ of cardinality at least $t'$. Observe that set $A$ and sets $A_e$, for each $e\in E(G)$, are minimal cut-modules and satisfy Lemma~\ref{lem:cuts}. We can then assume that all vertices in $A\cup\bigcup_{e\in E(G)} A_e$ are not in $X$. The number of the remaining vertices is then equal to 
\begin{align*}
|V(G')| - |A\cup\bigcup_{e\in E(G)} A_e| & = |V(G)|+|B|+|T|+|T_e|\cdot m \\
& = n+k+ (n+k)(m+1)=(n+k)(m+2). 
\end{align*}
Since $t'$ is at least $(m+1)(n+k)+2$, by the pigeonhole principle, we deduce that for each edge $e$, the pair of sets $T$ and $T_e$ are such that $T\cap X\not = \emptyset$ and $T_e\cap X\not = \emptyset$. We claim that $B\cap X=\emptyset$. If $|B\cap X|>1$, two vertices $u\in T$ and $v\in T_e$ cannot have $k+1$ disjoint shortest $u,v$-paths without internal vertices in $X$, for each $e\in E(G)$. Then $|B\cap X|\leq 1$. Since $t'\geq (m+1)(n+k)+2$ and $|T\cup\bigcup_{e\in E(G)} T_e|=(m+1)(n+k)$, there must be at least one vertex of $X$ in $V(G)$. Let $i\in V(G)\cap X$ be such a vertex and $e=ij$ be an edge incident to $i$. Under these conditions, a vertex $u\in T$ and a vertex $v\in T_e$ cannot be both in $X$, a contradiction. This concludes the proof that $B\cap X=\emptyset$. Since $t'=(m+1)(n+k)+t$, and $A$, $B$, and each set $T_e$, for each $e\in E(G)$, do not contain elements of $X$, then $t$ vertices of $X$ are necessarely in $V(G)$. These $t$ vertices form an independent set as, otherwise, two of them are incident to an edge $e$. But, in this case, again a vertex $u\in T$ and a vertex $v\in T_e$ cannot be both in $X$, a contradiction.

Finally, it is not difficult to check that $\diam(G') = 4$. 
\end{proof}

The above result immediately implies that the \ftmv problem is NP-hard.

\section{Graphs with large $\mv{k}$ numbers}
\label{sec:large}

In this section, we consider the graphs $G$ whose fault mutual-visibility number is large. More precisely, we describe the structure of the graphs $G$ with $\mv{k}(G) = n(G) - t$ for all $k\ge 1$ and $t\in \{0,1,\dots, k+1\}$. When $0\le t\le k$, the graphs with $\mv{k}(G) = n(G) - t$ coincide with the graphs whose clique number is $n(G) - t$. For the case $t = k+1$ a structural characterization is provided. 

\begin{theorem}
\label{thm:versus-omega}
If $G$ is a connected graph and $k\ge \ell \ge 0$, then $\mv{k}(G) = n(G) - \ell$ if and only if $\omega(G) = n(G) - \ell$.
\end{theorem}

\begin{proof}
    Assume first that $\mv{k}(G) = n(G) -\ell$. By \eqref{eq:omega-is-lower-bound} we have $\omega(G) \leq n(G) - \ell$. Suppose now on the contrary that $\omega(G) < n(G) - \ell$ and consider a $\mv{k}(G)$-set $X$. As $|X| = n(G)-\ell$ and $\omega(G) < n(G) - \ell$, there exist two vertices $u,v \in X$ with $uv \notin E(G)$. As there exist $k+1$ shortest $u,v$-paths without internal vertices in $X$, we get $\ell \ge k+1$. But this contradicts our assumption that $k\ge \ell$. 
    
    Conversely, assume that $\omega(G) = n(G) -\ell$. Then by \eqref{eq:omega-is-lower-bound} we know that $\mv{k}(G) \geq n(G) - \ell$. If $\mv{k}(G) > n(G) - \ell$ would hold, then by the above paragraph we would have $\omega(G) > n(G) - \ell$. We can conclude that $\mv{k}(G) = n(G) - \ell$.
\end{proof}

\begin{corollary}
\label{cor:n-and-minus-one}
If $G$ is a connected graph, the following properties hold. 
\begin{enumerate}
\item[(i)] If $k\ge 0$, then $\mv{k}(G) = n(G)$ if and only if $\omega(G) = n(G)$.
\item[(ii)] If $k\ge 1$, then $\mv{k}(G) = n(G) -1$ if and only if $\omega(G) = n(G) -1$.
\end{enumerate}
\end{corollary}

\begin{proof}
To get (i), set $\ell = 0$ in Theorem~\ref{thm:versus-omega}, and to get (ii) set $\ell = 1$. 
\end{proof}

With respect to Corollary~\ref{cor:n-and-minus-one}(ii) we add that a description of graphs with $\mu(G) = n(G) - 1$ is given in~\cite[Lemma 4.8]{distefano-2022}.

Theorem~\ref{thm:versus-omega} deals with situations in which $k\ge \ell$. The cases when $k < \ell$ are intrinsically different. For a small example consider stars $K_{1,m}$, for which we have $\mv{0}(K_{1,m}) = \mu(K_{1,m}) = m$. For any $k\ge 1$ we next characterize graphs $G$ with $\mv{k}(G) = n(G) - k - 1$ (that is, $\ell = k + 1$). To formulate the result we need the following family of graphs. 

\begin{definition}
\label{def:Hn}    
For $k\ge 1$ and $n\ge k+2$, the family of graphs $\mathcal{H}_{n,k}$ is defined as follows. Start with an arbitrary graph $H$ of order $n - k - 1$. Let $A$ be the set (possibly empty) of the universal vertices of $H$, and let $B = V(H)\setminus A$.  Next, let $W$ be a set of $k+1$ additional vertices, and add all possible edges between $W$ and $B$. Finally, add some arbitrary edges between the vertices from $W$, and some arbitrary edges between the vertices from $W$ and the vertices from $A$, where we take care that the obtained graph $\widehat{H}$ is connected and fulfils $\omega(\widehat{H})\le n-k-1$. The family $\mathcal{H}_{n,k}$ then contains all graphs $\widehat{H}$ constructed as described. 
\end{definition}

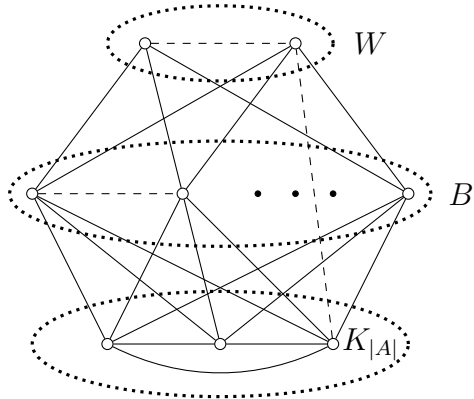
\begin{figure}[ht!]
	\begin{center}
		\begin{tikzpicture}[scale=1.0]
			\tikzstyle{vertex}=[circle, draw, inner sep=0pt, minimum size=6pt]
			\tikzset{vertexStyle/.append style={rectangle}}

                \vertex (1) at (0,0) [scale=0.7,fill=white] {};
                \vertex (2) at (2,0) [scale=0.7,fill=white] {};
                \vertex (3) at (5,0) [scale=0.7,fill=white] {};

                \node (B) at (5.7,0) {$B$};

                \vertex (3-1) at (3,0) [scale=0.35,fill=black] {};
                \vertex (3-2) at (3.5,0) [scale=0.35,fill=black] {};
                \vertex (3-3) at (4,0) [scale=0.35,fill=black] {};

                \vertex (4) at (1.5,2) [scale=0.7,fill=white] {};
                \vertex (5) at (3.5,2) [scale=0.7,fill=white] {};

                \vertex (6) at (1,-2) [scale=0.7,fill=white] {};
                \vertex (7) at (2.5,-2) [scale=0.7,fill=white] {};
                \vertex (8) at (4,-2) [scale=0.7,fill=white] {};

                \node (A) at (4.5,-2) {$K_{|A|}$};

                \path
                (1) edge (4)
                (1) edge (5)
                (2) edge (4)
                (2) edge (5)
                (3) edge (4)
                (3) edge (5)
                (6) edge (1)
                (6) edge (2)
                (6) edge (3)
                (7) edge (1)
                (7) edge (2)
                (7) edge (3)
                (8) edge (1)
                (8) edge (2)
                (8) edge (3)
                ;

                \path
                (8) edge (5)
                (1) edge (2)
                (4) edge (5)[dashed]
                ;

                \path
                (6) edge (7)
                (7) edge (8)
                ;

                \draw (6) .. controls (2,-2.5) and (3,-2.5) .. (8);

                \draw (2.5,-2) ellipse [x radius=2.5, y radius=0.7] [dotted, very thick];

                \draw (2.5,0) ellipse [x radius=2.8, y radius=0.7] [dotted, very thick];

                \draw (2.5,2) ellipse [x radius=1.5, y radius=0.5] [dotted, very thick];

                \node (W) at (4.5,2) {$W$};


	    \end{tikzpicture}
		\caption{The structure of the graphs from $\mathcal{H}_{n,1}$.}
		\label{fig:H-family}
	\end{center}
\end{figure}

For instance, selecting $H$ to be an edgeless graph of order $n-2$ we get $K_{2,n-2}\in \mathcal{H}_{n,1}$. A typical representative of the family $\mathcal{H}_{n,1}$ can be seen in Fig.~\ref{fig:H-family}. 

\begin{theorem}
\label{thm:n-2} 
If $k\ge 1$ and $G$ is a connected graph, then $\mv{k}(G) = n(G) - k - 1$ if and only if $G\in \mathcal{H}_{n(G),k}$.  
\end{theorem}

\begin{proof}
Let $G$ be a connected graph with $\mv{k}(G) = n(G) - k -1$. Let $X$ be a $\mv{k}(G)$-set of $G$, so that $|X| = n(G) - k - 1$, and let $V(G)\setminus X = \{w_1, \dots, w_{k+1}\} = W$. Let $A = \{u\in X:\ \deg_{G[X]}(u) = n(G) - k - 2\}$ and $B = X\setminus A$. If $u$ is an arbitrary vertex of $B$, then there exists a vertex $u'\in B$ such that $uu'\notin E(G)$. Since $X$ is a $k$-ftmv set, the only possibility to have $k+1$ internally disjoint shortest $u,v$-paths is that there exist the paths $u-w_i-u'$ for $i\in [k+1]$. It follows that every vertex from $B$ is adjacent to all the vertices of $W$. Note that a possible presence of some edges between the vertices of $W$ does not affect this conclusion. By~\eqref{eq:omega-is-lower-bound} we also have $\omega(G)\le n(G)-k - 1$. We can conclude that $G\in \mathcal{H}_{n(G),k}$.

Conversely, assume that $G\in \mathcal{H}_{n(G),k}$. If we set $G = \widehat{H}$, where $\widehat{H}$ is obtained from $H$ with the sets $A$ and $B$ as described in Definition~\ref{def:Hn}, then $A\cup B$ is a $k$-ftmv set of $G$, hence we have $\mv{k}(G) \ge  n(G) - k-1$. Suppose that $\mv{k}(G) \ge  n(G) - k$. Then $\mv{k}(G) =  n(G) - k+a$, where $0 \leq a \leq k$. In view of Theorem~\ref{thm:versus-omega} we can also write $\mv{k}(G) =  n(G) - \ell$, where $\ell = k-a$. Since clearly $k\ge k-a$, Theorem~\ref{thm:versus-omega} implies that $\omega(G)= n(G) - k+a \geq n(G) - k$, contradicting that $\omega(G)\le n(G)-k-1$ by to the definition of $G = \widehat{H}$. Hence, we can conclude that $\mv{k}(G) =  n(G) - k-1$.
\end{proof}

Related to Theorem~\ref{thm:n-2} we add that we are not aware of a characterization of graphs $G$ with $\mv{0}(G) = n(G) - 2$, that is, of graphs with the mutual visibility number equal to the order minus two. In view of~\eqref{eq:monotone}, this class forms a superclass of $\mathcal{H}_n$. The classes do not coincide as the graph $K_{3,3}-e$ demonstrates. Indeed, one can check that $\mv{0}(K_{3,3}-e) = 4$ and $\mv{1}(K_{3,3}-e) = 3$.

\section{Cartesian products}
\label{sec:cartesian}

In this section, we determine exact formulae for $\mv{k}(G)$, when $G$ is the  Cartesian product of two paths (grid graphs), the  Cartesian product of a path and a cycle (cylinder graphs), and the  Cartesian product of two cycles (tori graphs). Throughout the section we set $V(P_s) = [s]$, $s\ge 1$, and  $V(C_n) = [n]$, $n\ge 3$, and compute the indices modulo $n$. 

\subsection{Grids}

In this subsection, we consider the grids $P_m\cp P_n$, where we assume without loss of generality that $1\le m\le n$. 

$P_1 \cp P_1 \cong K_1$ for which $\mv{k}(K_1)=1$, $k\ge 0$. If $m = 1$ and $n \ge 2$, then $P_1 \cp P_n \cong P_n$, then we have $\mv{k}(P_n)=2$, $k \ge 0$. All the remaining cases are covered by the next result. 

\begin{theorem}\label{thm:grids}
If $2\le m\le n$, then
$$\mv{k}(P_m \cp P_n)=\left\{\begin{array}{ll}
                      2m; & k=0, \\
                      m; & k=1, \\
                      2; & k\ge 2.
                    \end{array}
\right.$$
\end{theorem}

\begin{proof}
In~\cite[Theorem 4.6]{distefano-2022} it has been shown that $\mu(P_m \cp P_n)= 2\min(m,n) = 2m$. Hence the result holds for $k=0$. Assume in the rest that $k\ge 1$. 

\medskip\noindent
{\bf Case 1}: $m=2$.\\
Let $n\ge 2$ and set $G=P_2\cp P_n$. We need to show that $\mv{k}(G)=2$. By~\eqref{eq:omega-is-lower-bound}, $\mv{k}(G)\ge 2$. Suppose now that $\mv{k}(G)\ge 3$. Then there are at least three vertices $u,v,w$ in a $\mv{k}$-set $X$ of $G$ and at least two of them, say $u,v$, are in the same $P_n$-layer. These two vertices must be adjacent since each layer is a convex subgraph isomorphic to $P_n$. Then $w$ is in the other layer. We may assume that $d(u,w)>d(v,w)$. Let $H$ the subgrid of $G$ having $u$ and $w$ on two opposite corners. Since $v$ is in $H$ and is adjacent to $u$, there is at most one shortest $u,w$-path without vertices in $X$, a contradiction. 

\medskip\noindent
{\bf Case 2}: $m=3$.\\
Let $G = P_3\cp P_n$, where $n\ge 3$. First, consider the subcase $k=1$. To see that $\mv{1}(G)\geq 3$, consider the set of vertices $\{(1,1), (2,2), (3,3)\}$. For the reverse inequality, suppose on the contrary that $\mv{1}(G)> 3$ and let $X$ be a $\mv{1}(G)$-set. Then at least one $P_n$-layer contains exactly two adjacent vertices $u,v$ from $X$. Now consider a third vertex $w\in X$ lying in another layer. Then, $\hull(\{u,v,w\})$ is a grid graph $H$ with a vertex $a\in\{u,v\}$ and $w$ in two corners of $H$. Since we cannot find two internally disjoint shortest $a,w$-paths without further vertices in $X$, we reach a contradiction. 

In the second subcase, let $k\ge 2$. By~\eqref{eq:omega-is-lower-bound}, $\mv{k}(G)\ge 2$. Suppose there exists a $k$-ftmv set $X$ with $|X|\ge 3$. Then $X$ contains two vertices $u,v\in X$ at distance at least two. Consider the subgrid $H$ with two corners $u$ and $v$. Since $H$ is a convex subgraph, each shortest $u,v$-path of $G$ belongs to $H$. Thus, since $\deg_{H}(u)= \deg_{H}(v)=2$, there cannot be more than two internally disjoint shortest $u,v$-paths, a contradiction. 

\medskip\noindent
{\bf Case 3}: $4\le m\le n$. \\
Consider first $k=1$. If $m$ is even, then consider the $m/2$ subgraphs induced by two consecutive $P_n$-layers isomorphic to $P_2\cp P_n$. These subgraphs are convex, hence by Case~1 and Lemma~\ref{lem:covering-by-convex-subgraphs} we have $\mv{1}(P_m \cp P_n) \leq m$. To show that $\mv{1}(P_m \cp P_n) \geq m$, consider the set of vertices $X = \{(i,i):\ i\in [m]\}$. If $m$ is odd, then consider the $(m-3) / 2$ subgraphs induced by two consecutive $P_n$-layers isomorphic to $P_2\cp P_n$, and one subgraph isomorphic to $P_3\cp P_n$. Then Cases 1 and 2, together with Lemma~\ref{lem:covering-by-convex-subgraphs} yield $\mv{1}(P_m \cp P_n) \leq m$, while the same set $X$ as in the even gives the reverse inequality. Assume second that $k \ge 2$. In this case, we proceed along the same lines as in the last paragraph of Case~2 to arrive at the conclusion that $\mv{k}(P_m \cp P_n)\le 2$. By~\eqref{eq:omega-is-lower-bound} we are done.
\end{proof}

\subsection{Cylinders}

In this section we consider $\mv{k}(P_m \cp C_n)$ for $k \ge 1$. The case $k=0$ has been considered earlier. In~\cite[Theorem 2.4]{korze-2024} it is proved that if $m+1\ge n\ge 6$, then $\mv{0}(P_m \cp C_n) = 2n$, and further conjectured in~\cite[Conjecture 2.6]{korze-2024} that if $n\ge 12$, then $\mv{0}(P_m \cp C_n) = \min\{3m, 2n\}$. 

We start by observing simple properties about cycle graphs $C_n$. We recall from~\cite[Lemma 2.8]{distefano-2022} that $\mv{0}(C_n)=3$. For $k\ge 1$, it is easy to note that 
$$
\mv{k}(C_n) = \left\{\begin{array}{ll}
       3; 	& n=3, \\                        
       2; 	& n\ge 4.
                    \end{array}
\right.
$$
In particular, for $n\ge 4$, let $X=\{u,v\}$ be a $\mv{k}(G)$-set of $C_n$. If $k\ge 2$, then $u$ and $v$ are adjacent in $C_n$ (cf. Proposition~\ref{prop:delta}). If $k=1$ and $n$ odd, then $u$ and $v$ are adjacent. If $k=1$ and $n$ is even, then $u$ and $v$ are either adjacent or antipodal in $C_n$.

\begin{lemma}\label{lem:prism-bounds}
If $m\ge 2$, $n\ge 4$, and $k\ge 1$ then
$$
\mv{k}(P_m \cp C_n) \le \left\{\begin{array}{ll}
       \min\{2m,n\}; 			& n \mbox{ even}, \\
       \max\{3,\min\{m,n\}\}; 	& n \mbox{ odd}.                      
                    \end{array}
\right.
$$
\end{lemma}
\begin{proof}
Observe first that each pair of adjacent $P_m$-layers induces a $P_m\cp P_2$ grid, and this grid is a convex subgraph of $P_m \cp C_n$. Since $\mv{k}(P_m \cp P_2) = 2$ (cf. Theorem~\ref{thm:grids}), Lemma~\ref{lem:covering-by-convex-subgraphs} immediately implies that $\mv{k}(P_m \cp C_n) \le n$ when $n$ is even. Assume next that $n$ is odd and let $X$ be a $\mv{k}(G)$-set of $P_m\cp C_n$. If each $P_m$-layer contains at most one vertex of $X$, then we are done. Assume thus that there exists $i\in [n]$ such that $|V(P_m^i)\cap X| = 2$. Then Lemma~\ref{lem:covering-by-convex-subgraphs} implies that $V(P_m^{i-1})\cap X = \emptyset$ and $V(P_m^{i+1})\cap X = \emptyset$. Then $n-3$ is even, and a repetitive application of Lemma~\ref{lem:covering-by-convex-subgraphs} implies that in these $P_m$-layers there are at most $n-3$ vertices from $X$, hence in this subcase $|X|\le n-1$. 

Assume $n$ even. Each $C_n$-layer is a convex subgraph of $P_m \cp C_n$. Since $\mv{k}(C_n) =2$, Lemma~\ref{lem:covering-by-convex-subgraphs} gives $\mv{k}(P_m \cp C_n) \le 2m$ when $n$ is even. By combining this bound with that in the first part of the proof, we get $\mv{k}(P_m \cp C_n) \le \min\{2m,n\}$.
 
Assume now that $n\ge 5$ odd. Consider first the case in which a $k$-ftmv set $X$ of $P_m \cp C_n$ contains two adjacent vertices $u=(i,j)$ and $v=(i,j+1)$ belonging to the same ${^i}C_n$-layer. Let $w=(i,j')$ be the vertex in the ${^i}C_n$-layer  at the same distance from both $u$ and $v$. If $w'\in X\setminus \{u,v\}$, then $k\ge 1$ implies that $w$ and $w'$ necessarily belong to the same $P_m^{j'}$-layer (i.e. $w'=(i',j')$, with $i'\neq i$). Assume that $X$ contains an additional fourth vertex $w''\in X\setminus \{u,v,w'\}$. As $w'$, also $w''$ must belong to $P_m^{j'}$-layer. Without loss of generality, assume $d(w'',v) > d(w',v)$. In this situation, it is easy to observe that from $w''$ to $v$ there exists only one shortest path with no internal vertices in $X$. We can conclude that, assuming $u,v\in X$ with $u$ and $v$ adjacent and belonging to the same $C_n$-layer, leads to $|X| \le 3$. Conversely, assuming that each vertex in $X$ belongs to a distinct $C_n$-layer, leads to $|X|\le m$. By combining these two bounds with that in the first part of the proof, and by using the construction in which $|X|=3$, we get the requested bound $\mv{k}(P_m \cp C_n) \le \max\{3,\min\{m,n\}\}$.
\end{proof}

We can now prove the following result about cylinders.

\begin{theorem} 
If $m\ge 2$ and $n\ge 3$, then
$$\mv{k}(P_m \cp C_n)=\left\{\begin{array}{ll}
    3; 						& k\ge 1, m\ge 2, n=3, \\                        
    \min\{2m,n\}; 			& k=1, m\ge 2, n\ge 4 \mbox{ even}, \\                  
    \max\{3,\min\{m,n\}\}; 	& k=1, m\ge 2, n\ge 5 \mbox{ odd}, \\     
    2; 						& k\ge 2, m\ge 2, n\ge 4.
    \end{array}
\right.$$
\end{theorem}

\begin{proof}

We distinguish several cases. 

\medskip\noindent 
{\bf Case $k\ge 1, m\ge 2, n=3$}.\\ 
Let $X$ be a $k$-ftmv set. Lemma~\ref{lem:prism-bounds} states that, regardless the parity of $n$, $|X|\le n$. Since $n=3$, the statement follows by letting $X$ contain the three vertices of any $C_n$-layer (i.e., a clique).

\medskip\noindent 
{\bf Case $k=1, m\ge 2, n\ge 4 \mbox{ even}$}. \\
From Lemma~\ref{lem:prism-bounds} we get 
$\mv{1}(P_m \cp C_n) \le \min\{2m,n\}$. When $n\le 2m$, set 
$$X= \{
(1,1), (1,\frac{n}{2}+1), 
(2,2), (2,\frac{n}{2}+2), 
\ldots, 
(n/2,n/2), (n/2,n) 
\}.$$ 
When $2m < n$, set 
$$X= \{
(1,1), (1,\frac{n}{2}+1), 
(2,2), (2,\frac{n}{2}+2), 
\ldots, 
(m,m), (m,\frac{n}{2}+m) 
\}.$$ 
Observe that $|X|=n$ in the first case, and $|X|=2m$ in the second one. The statement follows by easily checking that, in both cases, $X$ is a $1$-ftmv set ($n$ is even, in each $C_n$-layer there are two antipodal vertices, and in each $P_m$-layer there is exactly one vertex).

\medskip\noindent 
{\bf Case $k=1, m\ge 2, n\ge 5 \mbox{ odd}$}. \\
Analyze first the case $m=2$. From Lemma~\ref{lem:prism-bounds} we get 
$ \mv{1}(P_2 \cp C_n) \le 3$. In this case, the statement follows by taking $X=\{(1,1), (1,2), (2,\frac{n+1}{2})\}$ as a $1$-ftmv set.
In the remainder, assume $m\ge 3$. From Lemma~\ref{lem:prism-bounds} we get $ \mv{1}(P_m \cp C_n) \le \min\{m,n\}$. If $m\le n$, then set $X = \{(1,1), (2,2), \ldots, (m,m)\}$. If $m > n$, then set $X = \{(1,1), (2,2), \ldots, (n,n)\}$. In both cases, it is easy to observe that $X$ is a $1$-ftmv set of $P_m \cp C_n$ ($n$ is odd and in each layer there is exactly one vertex).

\medskip\noindent 
{\bf Case $k\ge 2, m\ge 2, n\ge 4$}.\\
For $k= 2$, suppose by contradiction there exists a $\mv{2}(P_m \cp C_n)$-set $X$ with three elements. There are two cases to be analyzed: (1) two elements of $X$ are in the same $C_n$-layer, and (2) the elements of $X$ are in distinct $C_n$-layers. In the first case, assume $u,v\in X$  are in the same ${^i}C_n$-layer and the third element of $X$ (say $w$) is in the ${^j}C_n$-layer, $j\neq i$. Since each $C_n$-layer is a convex subgraph and $k\ge 2 = \Delta(C_n)$, Proposition~\ref{prop:delta} applied to a ${^i}C_n$ implies that $u$ and $v$ must be adjacent. Since the degree of $u$ is three and one of its neighbors is in $X$, there are at most two internally disjoint $u,w$-shortest paths with no internal vertices in $X$. Hence, $X$ cannot be $2$-ftmv set of $G$. 
In the second case, assume $X=\{u,v,w\}$ and $u\in {^i}C_n$, $v\in {^j}C_n$, and $w\in {^\ell}C_n$, with $i<j<\ell$. If $n$ is odd, only one of the two neighbors $u'$ and $u''$ of $u$ in ${^i}C_n$ can belong to a shortest $u,v$-path. Conversely, if $n$ is even, both $u'$ and $u''$ belong to shortest $u,v$-paths, but such paths share the vertex that is antipodal to $u$ in ${^i}C_n$. Hence, $X$ cannot be a $2$-ftmv set of $G$. We can conclude that $\mv{2}(P_m \cp C_n) \le 2$ and thus $\mv{k}(P_m \cp C_n) \le 2$ for $k\ge 2$ by~\eqref{eq:monotone}. Therefore, $\mv{k}(P_m \cp C_n) = 2$ for $k\ge 2$. 
\end{proof}

\subsection{Tori}

In \cite[Corollary 4.7]{distefano-2022}, it is showed that $\mv{0}(C_m \cp C_n)\leq 3 \min\{m,n\}$ for any $m \geq 3$ and $n \geq 3$. In \cite[Theorem 3.8 and Table 2]{korze-2024}, Kor\v ze and Vesel determined the value of $\mv{0}(C_m \cp C_n)$ for all $m,n \geq 14$ and for all the other cases, that is, for $m \leq 13$ and $n \leq 14$, with computer assistance. See~\cite{korze-2024} for all the values. 

Since each $C_m$-layer and each $C_c$-layer is a convex subgraph of $C_m \cp C_n$, we can deduce from Lemma~\ref{lem:covering-by-convex-subgraphs} that for any $k \ge 1$, and any $m,n \ge 3$,
\begin{equation}
\label{eq:upper-bound-tori}
\mv{k}(C_m \cp C_n) \leq 2 \min\{m,n\}.
\end{equation}

\begin{theorem}
\label{thm:tori-C3-case}
    If $n\ge 3$, then
$$\mv{k}(C_3 \cp C_n)=\left\{\begin{array}{ll}
                      9; & k=0, n \geq 6,\\
                      7; & k=0, n \in \{4,5\}, \\
                      6; & k=0, n =3, \mbox{or } k=1, n \mbox{ even}, \\
                      4; & k=1, n \mbox{ odd}, \\
                      3; & k\ge 2.
                    \end{array}
\right.$$
\end{theorem}

\begin{proof}
Set $G =C_3 \cp C_n$. By \eqref{eq:omega-is-lower-bound}, we have $\mv{k}(G) \geq \omega(G)=3$ for any $k$. In \cite[Lemma 2.8]{distefano-2022}, it is showed that $\mv{0}(C_n)= 3$ for any $n \geq 3$. For the case $k = 1$, we recall that the $\mv{1}(G)$-sets of an odd cycle of order at least $5$ can only formed by two adjacent vertices, and the $\mv{1}(G)$-sets of an even cycle can only formed by two adjacent vertices or two antipodal vertices. For $k \geq 2$, the $\mv{k}(G)$-sets of a cycle of order at least $4$ can only be formed by two adjacent vertices.

\medskip\noindent 
{\bf Case $k=0$.} \\
By \cite[Table 2]{korze-2024}, we have $\mv{0}(C_3 \cp C_3) = 6$ and $\mv{0}(C_3 \cp C_n) = 7$ if $n \in \{4,5\}$. Assume that $n\ge 6$. By Lemma~\ref{lem:covering-by-convex-subgraphs}, we know $\mv{0}(C_3 \cp C_n) \leq 9$ for any $n$. On the other hand, for $i \in [3]$ we set
$$X=\left\{\begin{array}{ll}
                      \{(i,1),(i,\frac{n}{3}+1),(i,\frac{2n}{3}+1)\}; & 0\equiv n\ (\bmod\ {3}), \\
                      \{(i,1),(i,\frac{n-1}{3}+1),(i,\frac{2n-2}{3}+2)\}; & 1\equiv n\ (\bmod\ {3}), \\
                      \{(i,1),(i,\frac{n+1}{3}+1),(i,\frac{2n-1}{3}+1)\}; & 2\equiv n\ (\bmod\ {3}).
                    \end{array}
\right.$$
In any case, it can be checked that $X$ is a $0$-ftmv set by the definition, implying that $\mv{0}(G) \leq |X|=9$. Hence, we can conclude that $\mv{0}(G) =9$.

\medskip\noindent 
{\bf Case $k=1$ and $n=3$}. \\ 
Observe that $X=\{(1,3),(2,1),(2,2),(3,3)\}$ is a $1$-ftmv set of $G$. Thus, we have $\mv{1}(G) \geq |X|= 4$.

Conversely, let $X$ be a $\mv{1}(G)$-set. Since $|X|=\mv{1}(G) \ge 4$, there are at least two vertices of $X$ in some $C_3$-layer. Suppose that $X$ contains the three vertices of a $C_3$-layer. For each vertex $v \in V(G) \setminus X$, let $x \in X$ be a vertex in a different $C_3$-layer from $v$. We cannot find two shortest paths between $x$ and $v$ without internal vertices in $X$, so that $v$ cannot be added to $X$. In this case, we have $|X| \leq 3$, a contradiction.

It follows that each $C_3$-layer has at most two vertices of $X$, and there exists a $C_3$-layer containing exactly two vertices of $X$, without loss of generality, say $(1,3)$ and $(3,3)$. By the choice of $X$ and by symmetry, only the vertices $(2,2)$ and $(2,1)$ can be added to $X$. Thus, we have $|X| \leq 4$.

\medskip\noindent 
{\bf Case $k=1$ and $n\ge 5$ is odd}. \\
Observe that $X=\{(1,\frac{n+1}{2}),(2,1),(2,n),(3,\frac{n+1}{2})\}$ is a $1$-ftmv set of $G$. Thus, we have $\mv{1}(G) \geq |X|= 4$.

Conversely, let $X$ be a $\mv{1}(G)$-set. Since $|X|=\mv{1}(G) \ge 4$, there are at least two vertices of $X$ in some $C_n$-layer. Furthermore, since $C_n$ is an odd cycle, there are exactly two adjacent vertices of $X$ in some $C_n$-layer. Without loss of generality, say them $(2,1)$ and $(2,n)$. By the choice of $X$ and by symmetry, only the vertices $(1,\frac{n+1}{2})$ and $(3,\frac{n+1}{2})$ can be added to $X$. Thus, we have $|X| \leq 4$.

\medskip\noindent 
{\bf Case $k=1$ and $n$ is even}. \\
Clearly, $n \geq 4$. Let $X=\{(i,1),(i,\frac{n}{2}+1)~|~ i\in [3]\}$. Then $X$ is the set of vertices consisting of all the first and the $(\frac{n}{2}+1)$-th $C_3$-layer. Note that every $C_n$-layer is an even cycle, and by symmetry it only suffices to verify that there are two shortest paths between $(2,1)$ and $(1,\frac{n}{2}+1)$ without internal vertices in $X$. We remark that one can pass $(2,2),(1,2),\dots,(1,\frac{n}{2})$ and another can pass $(2,n),(1,n),\dots,(1,\frac{n}{2}+2)$. Combining with \eqref{eq:upper-bound-tori}, we have $6=|X| \leq \mv{1}(G) \leq 2 \times 3=6$.

\medskip\noindent 
{\bf Case $k \geq 2$}. \\ 
For the case $k\geq 4$, the result immediately follows from Proposition~\ref{prop:delta}. Let $X$ be a $\mv{k}(G)$-set. Consider the case $k = 3$. For each two non-adjacent vertices $x_1,x_2 \in X$, there are at most three internally disjoint shortest $x_1,x_2$-paths in $G$. Hence, $G[X]$ must be a clique, and we have $3 \leq \mv{3}(G)=|X| \leq 3$.

Consider $k=2$. If $n$ is odd, then for each two non-adjacent vertices $x_1,x_2 \in X$, there are at most two internally disjoint shortest $x_1,x_2$-paths in $G$. Thus, we have $\mv{2}(G)= 3$.

If $n$ is even, we suppose that $|X|=\mv{2}(G) \ge 4$. Then there exists a $C_n$-layer containing exactly two adjacent vertices in $X$. Without loss of generality, say them $(2,1)$ and $(2,n)$. For any vertex $v \in V(G) \setminus X$, there are no three shortest paths between $v$ and $(2,1)$ or $(2,n)$ without internal vertices in $X$. Thus, any vertex of $V(G) \setminus X$ cannot be added to $X$, and we will obtain a contradiction that $|X| \leq 2$.
\end{proof}

As already mentioned, the values from Theorem~\ref{thm:tori-C3-case} for the case $k=0$ were earlier determined in~\cite{korze-2024}. However, since this was done with the help of a computer, we have included here a more explicit description of maximum mutual-viability set. 

\begin{theorem}
\label{thm:tori}
If $n \geq m \ge 4$, then
    $$\mv{1}(C_m \cp C_n)=\left\{\begin{array}{ll}
                      6; & m = 4, n \mbox{ odd},  \\
                      m; & m \geq 5, n \mbox{ odd},  \\
                      n; & m \mbox{ odd}, n < 2m \mbox{ even},  \\
                      2m; & \mbox{otherwise}.
                    \end{array}
\right.$$
Moreover, if $k \geq 2$ and $n \geq m \ge 4$, then
$$\mv{k}(C_m \cp C_n)=\left\{\begin{array}{ll}
                      4; & k=2, m, n \mbox{ even},  \\
                      3; & k=2, m \mbox{ even}, n \mbox{ odd}, \mbox{ or } m \mbox{ odd}, n \mbox{ even}, \\
                      2; & k=2, m,n \mbox{ odd}, \\
                      2; & k \geq 3.
                    \end{array}
\right.$$
\end{theorem}

\begin{proof}
    Set $G=C_m \cp C_n$ and consider the following cases.

    \medskip\noindent 
    {\bf Case $k=1$, $m=4$, and $n$ is odd}. \\ 
    Let $X=\{(1,1), (1,2), (3,1), (3,2),(2,\frac{n+1}{2}+1), (4,\frac{n+1}{2}+1)\}$. Then $X$ is a $1$-ftmv set of $G$, implying that $\mv{1}(G) \geq |X| \geq 6$. Conversely, let $X$ be a $\mv{1}(G)$-set. Since $\mv{1}(G)=|X| \ge 6$ and $n$ is odd, there exist two vertices of $X$ that are adjacent in some $C_n$-layer. Without loss of generality, we may assume that they are $(1,1)$ and $(1,2)$. For the $^1C_n$-layer and the $^2C_n$-layer, only $(2,\frac{n+1}{2}+1)$ can be added to $X$. For the $^3C_n$-layer and the $^4C_n$-layer, they have at most four vertices in $X$. Similarly, if one of them has two adjacent vertices in $X$, then another layer only can have at most one vertex in $X$. Thus, we can conclude that $\mv{1}(G)=|X| \leq 6$. Both sides prove the argument.
    
    \medskip\noindent 
    {\bf Case $k=1$, $m=4$, and $n$ is even}. \\
    Let $X=\{(1,1), (2,1), (3,2), (4,2),(1,\frac{n}{2}+1), (2,\frac{n}{2}+1), (3,\frac{n}{2}+2), (4,\frac{n}{2}+2)\}$. Then $X$ is a $1$-ftmv set of $G$, implying that $\mv{1}(G) \geq |X| \geq 8$. Combining with \eqref{eq:upper-bound-tori}, we have $\mv{1}(G) \leq 2 \times 4=8$.

\medskip
    In the rest of the proof, assume that $n \geq m \geq 5$ when $k =1$. We first claim the following. 
    \begin{itemize}\label{claim:one-star}
    \item[$(\star)$] If $n$ is odd, and $X$ is a $1$-ftmv set of $G$ such that two vertices of $X$ are adjacent in some $C_n$-layer, then
    $$|X|\leq \left\{\begin{array}{ll}
                      4; & m \mbox{ odd}, \\
                      6; & m \mbox{ even}.
                    \end{array}
    \right.$$
    \end{itemize}
    For some $i \in [m]$ and $j \in [n]$, let $(i,j)\in X$ and $(i,j+1)\in X$. In the first step, we can verify that the other vertices of the $^iC_n$-layer cannot be added to $X$. In the second step, we can verify that if $m$ is odd, then the other vertices of the $C_m^j$-layer and the $C_m^{j+1}$-layer cannot be added to $X$, and if $m$ is even, then only $(i+\frac{m}{2},j)$ and $(i+\frac{m}{2},j+1)$ among the vertices of these two layers can be added to $X$. For the third, if $m$ is odd, then there are at most two vertices of the $C_m^{j-\frac{n-1}{2}}$-layer which can be be added to $X$, implying that $|X| \leq 4$. If $m$ is even, then only at most two vertices of the $C_m^{j-\frac{n-1}{2}}$-layer and at most two vertices of the $^{(i+\frac{m}{2})}C_n$-layer can be be added to $X$, implying that $|X| \leq 6$.

    For any $n \ge m \ge 5$, it is easy to check that $\{(i,i)~|~ i\in [m]\}$ is a $1$-ftmv set of $G$. Thus, we have $\mv{1}(G) \geq m \geq 5$. This argument together with \hyperref[claim:one-star]{($\star$)} provides the following:
    \begin{itemize}\label{claim:two-star}
    \item[$(\star\star)$] If $n$ is odd, $m \neq 6$, and $X$ is a $1$-ftmv set of $G$, then every $C_n$-layer has at most one vertex in $X$.
    \end{itemize}
    
    \noindent\noindent 
    {\bf Case $k=1$, $n \geq 2m$ is even}. \\
    Consider the set of vertices $X=\{(i,i),(i,\frac{n}{2}+i)~|~ i\in [m]\}$. For any two vertices $u$ and $v$ in $X$, consider a subgrid $H'$ with two corners $u$ and $v$. It is easy to verify that there are at least two shortest $u,v$-paths of $H'$ (also of $G$) without internal vertices in $X$. Thus, $X$ is a $1$-ftmv set, and we have $\mv{1}(G) \geq |X|= 2m$. Combining with the proved bound \eqref{eq:upper-bound-tori}, we have $\mv{1}(G) = 2m$.

    \medskip\noindent 
    {\bf Case $k=1$, $n \leq 2m-2$ is even, and $m$ is odd}. \\ 
    Consider the set of vertices $X=\{(i,i)~|~ i\in [m]\} \cup \{(m+1-\frac{n}{2},m+1), (m+2-\frac{n}{2},m+2), \dots, (\frac{n}{2},n)\}$. It is checked that there are two shortest $u,v$-paths of $G$ without internal vertices in $X$. By the definition, we know that $X$ is a $1$-ftmv set and
    \begin{equation}
    \label{eq:1-ftmv=n}   
    \mv{1}(G) \geq |X|= n.
    \end{equation}

    On the other hand, let $X$ be a $\mv{1}(G)$-set. If every $C_m$-layer has at most one vertex in $X$, then $\mv{1}(G) = |X| \leq n$. Assume that there are two vertices of $X$ in some $C_m$-layer. Since $m$ is odd, such a $C_m$-layer contains exactly two adjacent vertices of $X$. If $m=5$, then $n \in \{6,8\}$. Similar to \hyperref[claim:one-star]{($\star$)}, we know $\mv{1}(G)=|X| \leq 6$. Combining with \eqref{eq:1-ftmv=n}, in this case we have $n \neq 8$ and $\mv{1}(G)= 6=n$. If $m \geq 7$, then $n \geq 8$. Similar to \hyperref[claim:one-star]{($\star$)}, we know that $\mv{1}(G)=|X| \leq 6<n$, a contradiction to \eqref{eq:1-ftmv=n}.

    \medskip\noindent 
    {\bf Case $k=1$, $n \leq 2m-2$ is even, and $m$ is even}. \\ 
    Now $m \geq 6$. Consider the set of vertices $X=\{(i,i), (i,\frac{n}{2}+i)~|~ i\in [\frac{n}{2}]\} \cup \{(\frac{n}{2}+1,\frac{n-m}{2}+1), (\frac{n}{2}+1,n+1-\frac{m}{2}), (\frac{n}{2}+2,\frac{n-m}{2}+2), (\frac{n}{2}+2,n+2-\frac{m}{2}), \dots, (m,\frac{m}{2}), (m,\frac{n+m}{2})\}$. It can be checked that $X$ is a $1$-ftmv set and henceforth, $\mv{1}(G) \geq |X|= 2m$. Combining with \eqref{eq:upper-bound-tori}, we can conclude that $\mv{1}(G)=2m$.

    \medskip\noindent 
    {\bf Case $k=1$, $n$ is odd, and $m$ is even}. \\ 
    First, consider the set of vertices $X=\{(i,i)~|~ i\in [m]\}$. Then $X$ is a $1$-ftmv set of $G$, implying that $\mv{1}(G) \geq |X|=m$. On the other hand, let $X$ be a $\mv{1}(G)$-set. If every $C_n$-layer has at most one vertex in $X$, then $\mv{1}(G) = |X| \leq m$. Otherwise there are two vertices of $X$ in some $C_n$-layer. Since $n$ is odd, such a $C_n$-layer contains exactly the two vertices of $X$, and these two vertices are adjacent. If $m=6$, then by \hyperref[claim:one-star]{($\star$)}, we have $\mv{1}(G)=|X| \leq 6=m$; if $m \geq 8$, then by \hyperref[claim:one-star]{($\star\star$)}, we have $\mv{1}(G)=|X| \leq m$.

    \medskip\noindent 
    {\bf Case $k=1$, $n$ is odd, and $m$ is odd}. \\ 
    The set $\{(i,i)~|~ i\in [m]\}$ is a $1$-ftmv set of $G$, and we have $\mv{1}(G) \geq m$. On the other hand, if $X$ is a $\mv{1}(G)$-set, then by \hyperref[claim:one-star]{($\star\star$)} we also get  $\mv{1}(G)=|X| \leq m$.

    \medskip\noindent 
    {\bf Case $k=2$, $m$ and $n$ are even}. \\ 
    Note that $\{(1,1),(m,1),(\frac{m}{2},\frac{n}{2}+1),(\frac{m}{2}+1,\frac{n}{2}+1)\}$ is a $2$-ftmv set of $G$, hence $\mv{2}(G) \geq 4$. Conversely, let $X$ be a $\mv{2}(G)$-set. If $X$ is independent, then without loss of generality, we may assume that $(1,1) \in X$. Clearly, the vertices in the $C_m{^1}$-layer and the $^1C_n$-layer cannot be added to $X$. Since for each $x \in X$, $x\ne (1,1)$, there are at least three shortest paths of $G$ between $x$ and $(1,1)$, we infer that $x \in V(C_m^{\frac{n}{2}+1}) \cup V(^{\frac{m}{2}+1}C_n) \setminus \{(1,\frac{n}{2}+1),(\frac{m}{2}+1,1)\}$. Moreover, we can see that $X$ intersects $V(C_m^{\frac{n}{2}+1}) \cup V(^{\frac{m}{2}+1}C_n) \setminus \{(1,\frac{n}{2}+1),(\frac{m}{2}+1,1)\}$ in at most one vertex. This implies that $\mv{2}(G) \leq 2$, a contradiction to the lower bound.

    Assume second that $X$ contains two adjacent vertices. Without loss of generality, we may assume that the two vertices are $\{(\frac{m}{2},\frac{n}{2}+1),(\frac{m}{2}+1,\frac{n}{2}+1)\}$. Then the vertices in the $C_m^{\frac{m}{2}}$-layer, the $C_m^{\frac{m}{2}+1}$-layer, and the $^{\frac{n}{2}+1}C_n$-layer do not belong to $X$. For each $x \in X$, $x \in V(C_m^{1}) \setminus \{(\frac{m}{2},1),(\frac{m}{2}+1,1)\}$. Clearly, $X$ intersects $V(C_m^{1}) \setminus \{(\frac{m}{2},1),(\frac{m}{2}+1,1)$ in at most two vertices. This implies that $\mv{2}(G) \leq 4$, and we have done.

    \medskip\noindent 
    {\bf Case $k=2$, one of $m,n$ is odd and the other is even}. \\ 
    If $m$ is even and $n$ is odd, then $\{(1,1),(\frac{m}{2}+1,\frac{n+1}{2}),(\frac{m}{2}+1,\frac{n+3}{2})\}$ is a $2$-ftmv set of $G$, and we have $\mv{2}(G) \geq 3$. Similarly, if $n$ is odd and $n$ is even, then $\{(1,1),(\frac{m+1}{2},\frac{n}{2}+1),(\frac{m+3}{2},\frac{n}{2}+1)\}$ is a $2$-ftmv set of $G$, and we have $\mv{2}(G) \geq 3$.
    
    Conversely, let $X$ be a $\mv{2}(G)$-set. Based on the above analysis, we may assume that $m$ is even and $n$ is odd. If $X$ is independent and assuming without loss of generality that $(1,1) \in X$, then for each $x \in X$ we have $x \in V(^{\frac{m}{2}+1}C_n)$. Clearly, $|X \cap V(^{\frac{m}{2}+1}C_n)| \leq 2$. This implies that $\mv{2}(G) \leq 3$.

    Assume that two vertices of $X$ are adjacent in $G$. Without loss of generality, we may assume that the two vertices are $\{(\frac{m}{2}+1,\frac{n+1}{2}),(\frac{m}{2}+1,\frac{n+3}{2})\}$ or $\{(\frac{m}{2},\frac{n+1}{2}),(\frac{m}{2}+1,\frac{n+1}{2})\}$. In the first case, for each $x \in X$ we have $x \in V(^1C_n)$. The set $X$ intersects $V(^1C_n)$ in at most one vertex, implying that $\mv{2}(G) \leq 3$. In the second case, each vertex in $V(G) \setminus \{(\frac{m}{2},\frac{n+1}{2}),(\frac{m}{2}+1,\frac{n+1}{2})\}$ cannot be added to $X$, implying that $\mv{2}(G) \leq 2$, a contradiction.

    \medskip\noindent 
    {\bf Case $k=2$, $m,n$ are odd}. \\ 
    The lower bound $\mv{2}(G) \geq 2$ is obvious. Let $X$ be a $\mv{2}(G)$-set. Similar to the two cases above, assume that $(1,1) \in X$. Then each $x \in X$ must be a neighbor of $(1,1)$, hence $X$ intersects the open neighborhood of $(1,1)$ in at most one vertex. This implies that $\mv{2}(G) \leq 2$.

    \medskip\noindent {\bf Case $k\geq 3$}. \\ 
    The lower bound $\mv{k}(G) \geq 2$ is obvious. Let $X$ be a $\mv{k}(G)$-set. Assume that $(1,1) \in X$. For each $x\in X$, $x\ne (1,1)$, there are at least four shortest paths of $G$ between $x$ and $(1,1)$. Because of this we can check that when $m,n$ are even, only the vertex $(\frac{m}{2}+1,\frac{n}{2}+1)$ or only one neighbor of $(1,1)$ can belong to $X$ besides $(1,1)$. For other subcases, each $x \in X$ is a neighbor of $(1,1)$ in $G$, and $X$ intersects the open neighborhood of $(1,1)$ in at most one vertex. This implies that $\mv{k}(G) \leq 2$.
\end{proof}

\section{Graphs of diameter two}

In this section, we consider some classes of graphs of diameter two, including the Cartesian product and the direct product of complete graphs.

As a preliminary example, consider the Petersen graph $P$. Since any two vertices of $P$ at distance two are connected by a unique shortest path, we have $\mv{k}(P) = 2$ for any $k\ge 1$.

For a complete picture of the situation with Hamming graphs $K_m\cp K_n$, we recall the following. Let $z(m, n; 2, 2)$ denote the maximum number of $1$s in a $m\times n$ binary matrix, where the matrix contains no $2\times 2$ submatrix consisting of four $1$s. This is an instance of the notorious Zarankiewicz's problem, see~\cite{west-2021} for more information on it. In~\cite[Corollary~3.7]{cicerone-2023} it was proved that $\mu(K_m\cp K_n) = z(m,n;2,2)$, $m, n\ge 2$. We can now state: 

\begin{theorem}\label{thm:Hamming}
If $m\ge 2$ and $n\ge 2$, then
$$\mv{k}(K_m \cp K_n) = 
\left\{\begin{array}{ll}
z(m,n;2,2); & k=0, \\
m + n - 2; & k=1, \\
\max\{m,n\}; & k\ge 2.
\end{array}
\right.$$
\end{theorem}

\begin{proof}
By the discussion before the proof, $\mv{0}(K_m \cp K_n) = \mu(K_m \cp K_n) = 
z(m,n;2,2)$. 

Consider next the case $k=1$, and let $V(K_\ell) = [\ell]$, so that $V(K_m \cp K_n) = [m] \times [n]$. Let $X$ be a $\mv{1}$-set of $K_m\cp K_n$. Then we infer that for every $i,i'\in [m]$, $i\ne i'$, and for every $j,j'\in [n]$, $j\ne j'$, we have \begin{equation}
\label{eq:squares-at-most-2}
|X\cap \{(i,j), (i',j), (i,j'), (i',j')\}| \le 2.
\end{equation} 
Consider the bipartite graph $G_X$ with the bipartition $[m] \cup [n]$, where $ij\in E(G_X)$ if and only if $(i,j)\in X$. Then~\eqref{eq:squares-at-most-2} implies that $G_4$ does not include any $P_4$ as a subgraph. We hence infer that $\mv{1}(K_m \cp K_n)$ is equal to the largest number of edges in a bipartite graph, without $P_4$ subgraphs, whose bipartition sets are of respective cardinalities $m$ and $n$. Such a bipartite graph is the disjoint union of stars. If two stars have their central vertices, say $i$ and $i'$, in the same bipartition set, then we can replace the two stars with a single star having the center $i$ and as the set of leaves the union of the leaves of the two original stars. In this way, the vertex $i'$ now has no neighbor and can be added as a leaf to some star having the center in the other bipartition set. This argument implies that the largest number of edges will be achieved if and only if $G_X$ consists of exactly two disjoint stars, having their centers in different bipartition sets. Such two stars contain $m+n-2$ vertices, proving the expression for the case $k=1$.  

Consider finally the cases $k\ge 2$. Since any two nonadjacent vertices of $K_m\cp K_n$ are connected by exactly two shortest paths, a $k$-ftmv set of $K_m\cp K_n$ contains no non-adjacent vertices. In view of~\eqref{eq:omega-is-lower-bound} we can thus conclude that $\mv{k}(K_m \cp K_n) = \omega(K_m \cp K_n) = \max\{m,n\}$.
\end{proof}

We have seen that for any $k\ge 1$ we have $\mv{k}(P) = 2$, and that $\mv{k}(K_m \cp K_n)$ is linear in $n$ and $m$. All these graphs are of diameter two. Our next result demonstrates that among such graphs, $\mv{k}(G)$ can also be the of the same order as the order of the considered graph $G$. The case $k=0$ of the result has been earlier proved in~\cite[Theorem 2.2]{cicerone-2024b}. 

\begin{theorem}\label{thm:direct-Hamming}
If  $k\ge 0$ and $n,m \ge k+5$, then
$$\mv{k}(K_m \times K_n) = mn - 4 - k\,.$$
\end{theorem}

\begin{proof}
Let $k\ge 0$ and $n,m \ge k+5$. 

To prove that $\mv{1}(K_m \times K_n) \le mn - 4 - k$, suppose on the contrary that there exists a $k$-fault-tolerant mutual-visibility set $X$ of $K_m \times K_n$ of cardinality $mn - 3 - k$. Let $\overline{X} = V(K_m \times K_n) \setminus X$. Then $|\overline{X}| = 3 + k$. 

Assume first that $\overline{X}$ contains vertices  $(i_1,j_1), (i_2,j_2), (i_3, j_3)$, such that $|\{i_1, i_2, i_3\}| = 3$ and $|\{j_1, j_2, j_3\}| = 3$. Then we claim that at least one of the vertices $(i_2,j_1)$ and $(i_3,j_1)$ belongs to $X$. Since $m\ge k+5$ and $|\overline{X}| = 3 + k$, there exists two vertices $x$ and $x'$ from $V(K_m^{j_1})\cap X$. If all three vertices $(i_1,j_1)$, $(i_2,j_1)$, and $(i_3,j_1)$ would belong to $\overline{X}$, then $x$ and $x'$ would have at most $k$ common neighbors in $\overline{X}$, a contradiction. So assume without loss of generality that $(i_2,j_1)\in X$. If $(i_3,j_1)\in \overline{X}$, then $(i_2,j_1)$ and an arbitrary other vertex from $V(K_m^{j_1})\cap X$ again have at most $k$ common neighbors in $\overline{X}$ which is not possible. Hence we also have $(i_3,j_1)\in X$. But now $(i_2,j_1)$ and $(i_3,j_1)$ have have at most $k$ common neighbors in $\overline{X}$. This final contradiction proves that the case when $|\{i_1, i_2, i_3\}| = 3$ and $|\{j_1, j_2, j_3\}| = 3$ is not possible. 

If the vertices of $\overline{X}$ are contained in a single $K_m$-layer (or in a single $K_n$-layer for that matter), then two vertices from such a layer which are from $X$, have no common neighbors in $\overline{X}$ which is clearly not possible. So we are left with the case when there exist different indices $j$ and $j'$ such that $\overline{X}\subseteq V(K_m^{j}) \cup V(K_m^{j'})$. As $|\overline{X}| = k + 3$ and $m\ge k+5$, there exists an index $i$ such that $(i,j)\in X$ and $(i,j')\in X$. But now these two (non-adjacent) vertices have no common neighbor in $\overline{X}$, the final contradiction. 

We have thus proved that $\mv{1}(K_m \times K_n) \le mn - 4 - k$. To prove the reverse inequality, consider the sets $\overline{A} = \{(i,i):\ i\in [k+4]\}$ and $A = V(K_m \times K_n) \setminus \overline{A}$. We claim that $A$ is a $k$-ftmv set of $K_m \times K_n$. Let's consider two arbitrary, non-adjacent vertices of $A$. We may assume without loss of generality that these vertices are $(i,j)$ and $(i',j)$. If $j > s$, then each of $(i,j)$ and $(i',j)$ has at least $k+3$ common neighbors in $\overline{A}$. This in turn implies that $(i,j)$ and $(i',j)$ have at least $k+2$ common neighbors in $\overline{A}$. In the second case assume that $j\in [k+4]$.  Now each of $(i,j)$ and $(i',j)$ has at least $k+2$ common neighbors in $\overline{A}$. However, since $(j,j)$ is a common non-neighbor of $(i,j)$ and $(i',j)$, these two vertices still have at least $k+1$ common neighbors in $\overline{A}$. We can conclude that $A$ is a $k$-ftmv set of $K_m\cp K_n$. 
\end{proof}

Theorem~\ref{thm:direct-Hamming} has been in the special case $k=0$ earlier proved in~\cite[Theorem 2.2]{cicerone-2024b}.

\section{Concluding remarks}

We conclude our exposition by pointing out some open questions that may be of interest for continuing this line of research.
\begin{itemize}
    \item 
    In~\cite{bilo-2025}, strong inapproximability results for computing $\mv{0}(G)$ are established. It would be interesting to determine whether these inapproximability bounds extend to $\mv{k}(G)$ for $k \ge 1$.
    \item 
    In this work, we studied fault-tolerant mutual visibility for a number of well-known graph topologies obtained through different variants of graph products. For these graphs, our results extend existing knowledge on $\mv{0}(G)$. Since a substantial literature already exists on the computation and exact values of $\mv{0}(G)$ for many special graph classes, a natural direction is to investigate $\mv{k}(G)$, $k \ge 1$, on those same classes.
    \item 
    It would also be interesting to provide a structural characterization of graphs $G$ satisfying $\mv{1}(G) = 2$. Note, however, that this class of graphs is fairly large. By Proposition~\ref{prop:cut-vertex}, the problem can be restricted to biconnected graphs; within this category one already finds cycles of length at least $4$, pairs of cycles sharing exactly one edge, grids $P_2 \cp P_m$, and the Petersen graph.
    \item 
    In view of Theorem~\ref{thm:direct-Hamming}, it is a natural problem to determine the value of $\mv{k}(K_m\times K_n)$ when one of $m$ and $n$ is at most $k+4$.
    Using a computer, we have obtained the $k$-fault-tolerant mutual-visibility number for $K_m\times K_n$, for all $m, n\in \{3,4,5\}$ and all relevant $k$ in each of the cases. The results are collected in Table~\ref{tab:small-direct}.
\end{itemize}

\begin{table}[ht!]
    \centering
    \begin{tabular}{|c||c|c|c|} 
    \hline
        $k$ & 0 & 1 & $\ge 2$ \\
        \hline
        $\mv{k}(K_3\times K_3)$ & 6 & 4 & 3 \\
    \hline
    \end{tabular} \\ \vspace*{5mm}
    \begin{tabular}{|c||c|c|c|c|c|} 
    \hline
        $k$ & 0 & 1 & 2 & 3 & $\ge 4$ \\
        \hline
        $\mv{k}(K_3\times K_4)$ & 9 & 8 & 6 & 4 & 3  \\
    \hline
    \end{tabular} \\ \vspace*{5mm}
    \begin{tabular}{|c||c|c|c|c|c|c|} 
    \hline
        $k$ & 0 & 1 & 2 & 3 & 4 & $\ge 5$ \\
        \hline
        $\mv{k}(K_4\times K_4)$ & 12 & 12 & 10 & 9 & 6 & 4  \\
    \hline
    \end{tabular} \\ \vspace*{5mm}
    \begin{tabular}{|c||c|c|c|c|c|c|c|} 
    \hline
        $k$ & 0 & 1 & 2 & 3 & 4 & 5 & $\ge 6$ \\
        \hline
        $\mv{k}(K_3\times K_5)$ & 12 & 10 & 10 & 8 & 5 & 5 & 3  \\
    \hline
    \end{tabular} \\ \vspace*{5mm}
    \begin{tabular}{|c||c|c|c|c|c|c|c|c|c|c|} 
    \hline
        $k$ & 0 & 1 & 2 & 3 & 4 & 5 & 6 & 7 & 8 & $\ge 9$ \\
        \hline
        $\mv{k}(K_4\times K_5)$ & 16 & 16 & 15 & 13 & 12 & 10 & 8 & 6 & 5 & 4 \\
    \hline
    \end{tabular} \\ \vspace*{5mm}
    \begin{tabular}{|c||c|c|c|c|c|c|c|c|c|c|c|c|} 
    \hline
        $k$ & 0 & 1 & 2 & 3 & 4 & 5 & 6 & 7 & 8 & 9 & 10 & $\ge 11$ \\
        \hline
        $\mv{k}(K_5\times K_5)$ & 21 & 20 & 20 & 17 & 17 & 16 & 13 & 12 & 10 & 9 & 6 & 5 \\
    \hline
    \end{tabular} 
   \caption{The values of $\mv{k}(K_m\times K_n)$ for $m, n\in \{3,4,5\}$ and all relevant $k$.}
    \label{tab:small-direct}
\end{table}

\section*{Acknowledgments}

This work was supported by the Slovenian Research and Innovation Agency (ARIS) under the grants  P1-0297, N1-0285, N1-0355, N1-0431.
%
It was also partially supported by the European Union - NextGenerationEU under the Italian Ministry of University and Research (MUR) National Innovation Ecosystem grant ECS00000041 - VITALITY - CUP J97G22000170005, by the Italian National Group for Scientific Computation (GNCS-INdAM), and by the Academic Research Project MONET (University of L'Aquila).

\section*{Declaration of interests}

The authors declare that they have no known competing financial interests or personal relationships that could have appeared to influence the work reported in this paper.

\section*{Data availability}

Our manuscript has no associated data.



\end{document}